\documentclass[12pt]{article}
\usepackage{eucal,  amssymb,amsmath,graphicx, amsfonts, latexsym}
\usepackage[T1]{fontenc}
\newtheorem{theorem}{Theorem}[section]

\begin{document} \parskip=5pt plus1pt minus1pt \parindent=0pt
\title{Stochastic epidemics in growing populations}
\author{Tom Britton and Pieter Trapman}
\date{\today}
\maketitle

\begin{abstract} 
Consider a uniformly mixing population which grows as a super-critical linear birth and death process. At some time an infectious disease (of SIR or SEIR type) is introduced by one individual being infected from outside. It is shown that three different scenarios may occur: 1) an epidemic never takes off, 2) an epidemic gets going and grows but at a slower rate than the community thus still being negligible in terms of population fractions, or 3) an epidemic takes off and grows quicker than the community eventually  leading to an endemic equilibrium.  Depending on the parameter values, either scenario 1 is the only possibility, both scenario 1 and 2 are possible, or scenario 1 and 3 are possible.  
 \end{abstract}

\section{Introduction}\label{sec-intro}

Consider a population in which an infectious disease is introduced. The question of what might happen is the topic of many scientific papers, relying on mathematical or statistical models, or simulations, and often also fitted to specific populations and/or diseases.

In the current paper we use stochastic models to answer this type of questions for a dynamic population model. We explicitly include that also the population changes randomly over time: old individuals die and new are born, in such a way that the population, on average, grows. Our population is modeled by a super-critical linear birth and death process. An infectious disease is then introduced into this population. The disease is of SIR (susceptible-infectious-removed) type (e.g.\ \cite{DHB12}) meaning that individuals are at first Susceptible. If they get infected they become Infectious, and after some time in the infectious stage they Recover and become immune and stay so forever. In Section \ref{ext} we extend the result to an SEIR epidemic where individuals upon infection at first become latent (Exposed but not yet infectious) for some time before they enter the Infectious state.

We study the simplified model where individuals mix (and hence also infect) uniformly in the community. We prove that there are three qualitatively different scenarios that may happen. Either the disease outbreak dies out quickly. If not, there may be a big outbreak infecting more and more people, but still infecting only a negligible minority of the (growing) population. The third scenario is where the epidemic outbreak grows faster than the population, leading to an endemic situation in which the fraction of infectives fluctuates around some fixed level: the endemic level. Depending on the parameter values either only the first scenario is possible, the first and second, or the first and third scenario.

Modeling of infectious diseases in growing populations have been modeled by several authors. Thieme (1992) \cite{Thie92}  and Li et al.\cite{Li99} study fatal deterministic epidemic models allowing the population to be growing. Diekmann et al. \cite[Section 4.4]{DHB12} study a model very similar to the current SIR model, also allowing for disease fatalities, the difference being that we here look at a stochastic model for a \emph{finite} population. Breban et al.\cite{BVB05} considers approximations of epidemics taking place on networks in growing populations.

The rest of the paper is structured as follows. In Section 2 the growing population model and the SIR epidemic model are  defined. In Section 3 we state the results for the model and give the proofs. In Section 4 we discuss how to extend the results to an SEIR epidemic model, and discuss interesting future problems.

\section{The models}\label{sec-models}

\subsection{The dynamic population model}\label{sec-pop}

The population model is defined to be a super-critical Markovian branching process, or more specific a super-critical linear birth-and-death (B-D) process. More precisely, the life lengths of individuals are exponentially distributed with death rate $\mu$, and during the life each individual gives birth to new individuals according to a Poisson process with intensity $\lambda$, all random quantities mentioned being independent of each other. We assume that $\lambda>\mu$. i.e.\ that the process is super-critical. An alternative description of the process, the B-D description, is by the rates at which the process increases by one and decreases by one. If $N(t)$ denotes the number of individuals at $t$, then the rate of increase is $\lambda N(t)$ and the rate of decrease is $\mu N(t)$. The fact that $N(t)$ appears linearly in both rates makes it a linear B-D process, and the fact that the net increase/decrease is $\lambda-\mu>0$ (by assumption) makes the process super-critical. The starting state is $N(0)=n$ for some strictly positive integer $n$.

It is well-known (e.g.\ \cite{Hacc05,Jage75}) that $N(t)$ has positive probability of growing to infinity, and if it does, it grows at the rate $N(t)\sim e^{(\lambda -\mu)t}$. By which we mean that if $n=1$ than  $N(t)/e^{(\lambda -\mu)t}$ converges to a random variable $X$ with an atom of probability $\mu/\lambda$ at $X=0$, and conditioned on not being 0, $X$ is exponentially distributed with intensity parameter $1-\mu/\lambda$. In what follows we are only interested in the part of the sample space where this happens. For that reason we assume $N(0)=n$ is large (whenever $n$ is important notations are equipped with an $n$-index). For any fixed $\varepsilon$ we can in fact choose $n_0$ large enough such that for any $n\ge n_0$ we have

1) $P(\lim_t N_n(t)=\infty)>1-\varepsilon$, and 

2) $P((1-\varepsilon)ne^{(\lambda-\mu)t}\le N_n(t)\le (1+\varepsilon)ne^{(\lambda-\mu)t},\ \forall t)>1-\varepsilon$.

Further, it is also known \cite[Corollary 6.8.2]{Jage75} that the age distribution among living individuals converges to an exponential distributed random variable with mean $1/\lambda$.  More precisely, the age $A$ of a randomly selected individual converges (as $t\to\infty$) to  an exponential distributed random variable with mean $1/\lambda$, so the death rate does not affect the asymptotic age distribution.

\subsection{The uniform Markovian dynamic SIR epidemic model}\label{sec-hom-epid}

We now define a uniformly mixing Markovian epidemic model ''on'' the dynamic population model described in the previous section. We assume an SIR (Susceptible-Infectious-Removed) epidemic and extend this to an SEIR epidemic model, also having an Exposed/latent state when infected and prior to the infectious state, in Section \ref{ext}, cf.\ \cite{DHB12}.

Initially ($t=0$) the population consists of $N(0)=n$ individuals. At this time an SIR infectious disease is introduced by infecting one uniformly chosen individual. The epidemic is Markovian: An infectious individual remains infectious for an exponential time with rate $\delta$ (unless of course it dies before). When the infectious period stops, he/she  recovers and becomes immune for the rest of his/her life. During the infectious period the individual has infectious contacts randomly in time according to a homogeneous Poisson process with rate $\gamma$. Each time this happens the contacted person is randomly selected among all living individuals. If the contacted person is susceptible he/she becomes infected and infectious -- otherwise the contact has no effect. There is no vertical transmission (i.e., all new born individuals are susceptible) and the disease has no effect on the death rate which hence still equals $\mu$.

Let $Z(t)=(S(t), I(t), R(t))$ respectively, denote the number of susceptible, infectious and recovered individuals at time $t$. The epidemic is initiated at $Z(0)=(S(0), I(0), R(0))=(n-1,1,0)$ (the index $n$ is added when the initial population size is important).
The possible events and their rates, when currently in state $Z(t)=(S(t),I(t),R(t))=(s,i,r)=z$ are given in Table \ref{Tab-jumps}.
\begin{table}\centering
\begin{tabular}{|l|c|c|c|}
  \hline
Event & State change $\ell$& Rate name & Rate\\
  \hline
  Birth & $(1,0,0)$ & $\beta_{(1,0,0)}(s,i,r)$ & $\lambda (s+i+r)$\\
Death of susceptible & $(-1,0,0)$ &  $\beta_{(-1,0,0)}(s,i,r)$ & $\mu s$\\
Death of infective & $(0,-1,0)$ & $\beta_{(0,-1,0)}(s,i,r)$ & $\mu i$\\
Death of recovered & $(0,0,-1)$ & $\beta_{(0,0,-1)}(s,i,r)$ & $\mu r$\\
Infection & $(-1,1,0)$ & $\beta_{(-1,1,0)}(s,i,r)$ & $\gamma i\cdot s/(s+i+r)$\\
Recovery & $(0,-1,1)$ & $\beta_{(0,-1,1)}(s,i,r)$ & $\delta i$\\
  \hline
\end{tabular}
 \caption{The uniform Markovian dynamic SIR epidemic: type of events, their state change $\ell$ (the old state $z=(s,i,r)$ is hence changed to $z-\ell$) and their rates.}\label{Tab-jumps}
\end{table}

To summarize, the population model has two parameters, the birth rate $\lambda$ and the death rate $\mu$ (where $\lambda>\mu$ is assumed), and the epidemic  has two parameters: the infection rate $\gamma$ and the recovery rate $\delta$.

\section{Results}\label{sec-results}
We divide the results into the initial phase of an epidemic, where the fraction infected is small, and those where the epidemic grows to an endemic state.
  
\subsection{Initial phase of the epidemic}\label{sec-initial}

We now study the epidemic during the initial phase (when the \emph{fraction} infected is still small). If we focus on the number of infectives $I(t)$ we see from the table above that $I(t)$ increases at rate $\gamma I(t)\cdot S(t)/N(t)$ (due to infection) and decreases at rate $(\delta+\mu)I(t)$ (due to recovery or death). As a consequence, when the fraction infected is still small, or equivalently the fraction susceptible $S(t)/N(t)$ is close to 1, then the number of infectives $I(t)$ behaves almost like a linear branching process with birth rate $\gamma I(t)$ and death rate $(\delta+\mu)I(t)$. In what follows we will make this approximation precise.

Since $S(t)/N(t)\le 1$ we can dominate the number of infectives by an ''upper'' linear branching process $I^{(U)}(t)$ with $I^{(U)}(0)=1$, and having linear birth rate $\gamma I^{(U)}(t)$ and linear death rate $(\delta + \mu) I^{(U)}(t)$. Since the birth rate for $I^{(U)}(t)$ is always as large as that of $I(t)$ and the death rates are the same, the two processes can be coupled such that $I^{(U)}(t)\ge I(t)$ almost surely for all $t$ using methods very similar to those of Ball and Donnely \cite{Ball95}.

Obtaining a ''good'' lower bound is somewhat more involved. We do this for the process $I(t)$ up until the first time that $S(t)/N(t)\le 1-\varepsilon_0$, where $\varepsilon_0 >0$ can be chosen arbitrary small. This method is inspired by \cite{Durr07}. Before this time-point $I(t)$ can be bounded from below by the linear B-D process $I^{(L)}(t)$ having the same death rate $(\delta+\mu)I^{(L)}(t)$, but with birth rate $\gamma (1-\varepsilon_0)I^{(L)}(t)$. All three process can be defined on the same probability space (using the same random life-lengths) such that $I^{(L)}(t)\le I(t)\le I^{(U)}(t)$ for all $t$ at least up until the first time that $I(t)+R(t)\ge \varepsilon_0 N(t)$.

The number of infectives in the epidemic process is, at least until a fraction $\varepsilon_0$ of all living individual have been infected, sandwiched between the two linear B-D processes, both having individual death rate $\mu+\delta$, and individual birth rate $\gamma (1-\varepsilon_0)$ and $\gamma$ respectively. We use this to prove what will happen during the early stages of the epidemic in the theorem below. First we define two important quantities for the epidemic/branching processes. The first is the basic reproduction number, which is defined as the expected number of infectious contacts an infected individual has. For the current model this number equals
\begin{equation}
R_0=\frac{\gamma}{\delta + \mu}.\label{R_0}
\end{equation}
This follows since an individual has contacts at rate $\gamma$ during the infectious period which lasts for a mean duration of $1/(\delta+\mu)$; the individual stops being infectious due to recovery or death. The Malthusian parameter $\alpha$ is defined as the exponential growth/decay rate the epidemic/branching process has. It is given as the solution to 
$$
1 = \int_0^\infty e^{-\alpha t} c(t) dt,
$$
where $c(t)$ is the expected rate at which an individual has infectious contacts $t$ time units after it was infected \cite{Jage75}. For the current model, the contact rate is $\gamma$ and the contact is infectious if the individual is still infected, which it is with probability $e^{-(\delta+\mu)t}$. Putting this into the equation above and solving for $\alpha$ gives the Malthusian parameter being equal to
\begin{equation}
\alpha=\gamma - (\delta+\mu). \label{Malthus}
\end{equation}
From Equations (\ref{R_0}) and (\ref{Malthus}) it follows that the basic reproduction number $R_0$ exceeds the value of 1 if and only if $\alpha$ exceeds 0. The interpretation of this general rule is that a major outbreak can occur if and only if $R_0>1$ and this also implies that the epidemic will in expectation  \emph{grow} exponentially. We now state our main theorem concerning the early stage of an outbreak. What may happen depends on the sign of $\alpha=\gamma -(\delta +\mu)$ and on its relation to the growth rate $\lambda-\mu$ of the population.

\begin{theorem}\label{mainthm}
Consider the uniform epidemic model in the growing population defined above, and let $(S_n(0), I_n(0), R_n(0))=(n-1,1,0)$. As before, let $I^{(U)}(t)$ denote the  linear birth and death process with individual birth rate $\gamma$ and death rate $(\delta +\mu)$, with $I^{(U)}(0)=1$. Then, as $n\to\infty$, we have the following results:

i) The epidemic process $\{I_n(t);0\le t\le t_1\}$ converges weakly to $\{I^{(U)}(t);0\le t\le t_1\}$ on any finite interval $[0,t_1]$.

ii) If $\alpha<0$ then $I^{(U)}$ is sub-critical (critical if $\alpha= 0$). For any $n$, $I_n(t)\to 0$ as $t\to\infty$ with probability 1.

iii) If $0<\alpha<\lambda-\mu$, then $I^{(U)}$ is super-critical. Let $\pi_n:=P(\lim_t I_n(t)= 0)$. Then  $\pi_n\to (\delta+\mu)/\gamma$ as $n \to \infty$. With the remaining probability $1-\pi_n$, $I_n$ grows exponentially: for any $n$, $I_n(t)\sim e^{\alpha t}$, but $I_n(t)/N_n(t)\to 0$ as $t$ grows.

iv) If $\alpha\ge \lambda-\mu$, then $I^{(U)}$ is super-critical. Let $\pi_n:=P(\lim_t I_n(t)= 0)$. Then $\pi_n\to (\delta+\mu)/\gamma$. With the remaining probability $1-\pi_n$,   $I_n(t)\to \infty$ as $t \to \infty$.
\end{theorem}
\emph{Remark}.
What happens after the initial growth in case \emph{iv)} is treated in the next section.

\emph{Proof}
\emph{i)} The first statement follows from the fact shown above the theorem that $I_n(t)$ can be sandwiched between two linear birth and death processes having arbitrary close behavior on every finite time interval.

\emph{ii)}  It follows from standard birth-and-death process theory that $I^{(U)}$ is sub-critical and hence goes extinct with probability 1. The coupling $I_n(t)\le I^{(U)}(t)$ applies forever which then gives the result.

\emph{iii)} Since we may consider arbitrary large $n$ and hence arbitrary small $\varepsilon_0$, both bounding B-D-processes are super-critical. Their probabilities of extinction  (i.e.\ minor outbreak) are obtained by first deriving the distribution of the number of infectious contacts that one infected has, i.e.\ the offspring distribution, and then deriving the extinction probability of the corresponding branching process using theory for branching processes \cite{Jage75}. In the current situation we have exponential duration of the infectious period, the infectious period stops at rate $\delta+\mu$, and during this period an infective infects others at rate $\gamma$ and $\gamma(1-\varepsilon)$ respectively. As is well-known \cite{DHB12} the offspring distribution is then geometric and the minor outbreak probability equals $(\delta+\mu)/\gamma$ and $(\delta+\mu)/(\gamma(1-\varepsilon_0))$ respectively, but since $\varepsilon_0$ may be made arbitrary small the limiting probability of a minor outbreak equals $(\delta+\mu)/\gamma$. With the remaining probability the two B-D-processes grow exponentially at rate $\gamma-(\delta+\mu)$ and $\gamma(1-\varepsilon_0)-(\delta+\mu)$ respectively. Since the population grows at a higher exponential rate $\lambda-\mu$ and we start with arbitrary large population $N_n(0)=n$, the ratio $S_n(t)/N_n(t)=1-(I_n(t)+R_n(t))/N_n(t)$ will with arbitrary large probability never go below $1-\varepsilon_0$ (so the coupling never breaks down) and eventually $I_n(t)/N_n(t)\to 0$.

\emph{iv)} The initial statements are identical to those in \emph{iii)}. The only difference is that the population growth rate of $I_n(t)$ now is greater than that of $N_n(t)$ (at least if $\gamma> \delta + \lambda$), so eventually $S_n(t)/N_n(t)$ will be below $1-\varepsilon_0$, and the coupling breaks down. In the next section we show that in that case $I_n(t) \to \infty$ as $t \to \infty$. \hfill $\blacksquare$

Case iii) of the theorem might bring up the following question: Can we couple  $ I^{(U)}(t)$ and $I_n(t)$ such that  $ I^{(U)}(t)-I_n(t)=0$ hold forever with positive probability?   Using the theory as presented in \cite{Ball95}, we note that $ I^{(U)}(t)-I_n(t)=0$ forever with probability $\prod_{t \in \mathcal{T}}  S_n(t-)/N_n(t-)$, where $\mathcal{T}$ is the sets of times at which $I^{(U)}(t)$ increases by one. This probability is positive if and only if $$\sum_{t \in \mathcal{T}} (1- S_n(t-)/N_n(t-)) = \sum_{t \in \mathcal{T}} (I_n(t-) +R_n(t-))/N_n(t-)<\infty.$$ From the theory of branching processes and part iii) of the theorem we know that if the epidemic does not go extinct, then $e^{-(\lambda -\mu)t} N_n(t)$, $e^{-\alpha t} I_n(t)$ and $e^{-\alpha t} (I_n(t) +R_n(t))$ converges with probability 1, to positive random variables (say respectively $W_1$, $W_2$ and $W_3$. Heuristic arguments now give that for every $n$ and a large enough constant $c$,  $\sum_{t \in \mathcal{T}} (I_n(t-) +R_n(t-))/N_n(t-) \leq  c \int_0^{\infty} \gamma W_2e^{\alpha t} W_3e^{\alpha t}/(W_3 e^{(\lambda-\mu)t}.$
This sum is finite if $\alpha^2< \lambda-\mu$, which suggests that the event $ I^{(U)}(t)-I_n(t)=0$ forever has positive probability only in a part of the parameter domain satisfying  $0<\alpha<\lambda-\mu$. 

\subsection{Endemicity}\label{sec-end}

In this section we still assume that $\lambda-\mu>0$ (super-critical population process) and we only consider the case $\alpha> \lambda-\mu$ (the case $\alpha=\lambda-\mu$ is omitted), so that the epidemic, in case it takes off, grows at a higher rate than the population. We focus only on the situation where $I_n(t)$ and $I^{(U)}(t)$  tend to infinity; the situation where they go extinct being treated in the previous section.

As mentioned in Section \ref{sec-pop} the population size $N_n(t)$ will be close to the deterministic function $n(t)=ne^{(\lambda-\mu)t}$ in the sense that if $n=N_n(0)$ is chosen large enough, then $|N_n(t)/n(t)-1|<\varepsilon$ for all $t$ with probability arbitrary close to 1 (see Section 2). Since the population grows we introduce the ''proportion''-processes $\bar S_n(t)=S(t)/n(t)$, $\bar I_n(t)=I(t)/n(t)$ and $\bar R_n(t)=R(t)/n(t)$, and $\bar Z_n(t)=(\bar S_n(t), \bar I_n(t), \bar R_n(t))$.

Below we show that the vector valued stochastic process $\bar Z_n(t)$ converges to a deterministic vector-process $\bar z(t)$. In order to obtain a non-trivial deterministic process it is necessary to start the ''proportion process'' with a positive, albeit small, fraction infectives. For this reason we assume that $\bar Z_n(0)=(\bar S_n(0), \bar I_n(0), \bar R_n(0))=n^{-1}(S(0), I(0), R(0))=(1-\varepsilon_1-\varepsilon_2, \varepsilon_1, \varepsilon_2)$ for some small but otherwise arbitrary $\varepsilon_1,\varepsilon_2 >0$. This starting point is arbitrary chosen and we do not claim that it is exactly here that the stochastic epidemic process will pass through when the number of infectives grows in comparison with the population. In what follows we show that the starting point has hardly effect on the state of the process after a long time as long as $\bar S_n(0)$ is close to 1 and $\bar I_n(0)$ and $\bar R_n(0)$ are close to 0.

The proof that $\bar Z_n(\cdot )$ converges to $\bar z(\cdot )$ uses methods from \cite[Chapter 5]{EK05}, see also \cite[Chapter 5]{AB00} applying the theory to epidemic processes, with the difference that now a time-inhomogeneous normation is used. The epidemic process is defined in terms of (stochastic) rates (see Table \ref{Tab-jumps}). As in \cite{EK05} this implies that the epidemic process can equivalently be defined using Poisson processes. Let, for each of the six types of jumps (specified by $\ell$) defined in Table \ref{Tab-jumps}, $Y_\ell(\cdot)$ denote independent standard (intensity 1) Poisson processes, and let $\hat Y_\ell (\cdot)$ be the corresponding centered processes, so $\hat Y_\ell (t)=Y_\ell (t)-t$. We start by writing the original epidemic process $Z(t)=(S(t),I(t),R(t))$ using Poisson processes:
\begin{equation}
Z(t)= Z(0) + \sum_\ell \ell Y_\ell \left( \int_0^t\beta_\ell (Z(u))du \right).\label{epid-poisson}
\end{equation}
This follows from \cite{EK05} but can also be explained heuristically. The process starts in the correct point, and at any given time point $t$ the rate of the different jumps are as defined in Table \ref{Tab-jumps} and whenever a particular jump $\ell$ occurs the process $Z(t)=(S(t),I(t),R(t))$ changes with the vector $\ell$. The different Poisson processes were defined independently, the dependencies lie in for how long they are observed.

We now rewrite Equation (\ref{epid-poisson}) in terms of $\bar Z_n(t)$. First we observe that $\beta_{\ell}(Z(u))= n(u) \beta_{\ell}(\bar Z_n(u))$ for all $\ell$, except $\ell=(-1,+1,0)$ and $\ell=(+1,0,0)$ for which we have $\beta_{(-1,+1,0)}(Z(u))=(n^2(u)/N(u)) \beta_{(-1,+1,0)}(\bar Z_n(u))$ and $\beta_{(1,0,0)}(Z(u))=N(u) \beta_{(1,0,0)}(\bar Z_n(u))$. However, since $N(u)/n(u)$ is arbitrary close to 1 with arbitrary large probability, the same relation approximately holds also for these $\ell$ (the error term is neglected in what follows). Using this we obtain
\begin{equation}
\bar Z_n(t)= \frac{Z(0)}{n(t)} + \frac{1}{n(t)}\sum_\ell \ell \hat Y_\ell \left( \int_0^t n(u)\beta_\ell (\bar Z_n(u))du \right)+  \int_0^te^{-(\lambda-\mu)(t-u)}\sum_\ell \ell\beta_\ell (\bar Z_n(u))du.
\label{barepid-poisson}
\end{equation}

The middle term on the right will become small as $n$ grows (recall that $n(t)=ne^{(\lambda-\mu)t}$) because a centered Poisson process divided by something proportional to the mean of the non-centered process becomes negligible as the mean increases by the strong law of large numbers. This suggests that $\bar Z_n(t)$ converges to a deterministic vector process $\bar z(t)$ defined by
\begin{equation}
\bar z(t)=\bar z_0e^{-(\lambda-\mu)t}+  \int_0^te^{-(\lambda-\mu)(t-u)} \sum_\ell \ell\beta_\ell (\bar z(u))du,\label{bar-z}
\end{equation}
with $\bar z_0=( 1-\varepsilon_1-\varepsilon_2, \varepsilon_1, \varepsilon_2)$.

We are now ready to state our theorem for case iv in Theorem \ref{mainthm}. the situation where the epidemic initially grows at a higher rate than the population.
\begin{theorem}
Assume that $\alpha> \lambda-\mu$ and that the epidemic process starts in $Z_n(0)=(n(1-\varepsilon_1-\varepsilon_2), n\varepsilon_1, n\varepsilon_2)$ and consider the ''proportion process'' $\bar Z_n(t)=Z_n(t)/ne^{(\lambda-\mu)t}$. Then, as $n\to\infty$, $\bar Z_n(t) \to \bar z(t)$ in probability, uniformly on any finite interval, where $\bar z(t)$ was defined in (\ref{bar-z}).
\end{theorem}
\emph{Proof.}
The proof follows the same ideas as Theorem 5.2 in \cite{AB00} originating from \cite{EK05}, so we leave out some details. We also leave out the error term stemming from the approximation $N(u)/n(u)\approx 1$. Define $\bar\beta_\ell=\sup_{\bar z \in [0,1]^3} \beta_{\ell}(\bar z)$ which is finite. Further, let $M<\infty$ satisfy $|\sum_\ell \ell (\beta_\ell (x)-\beta_\ell (y))|\le M|x-y|$ which is true since all $\beta_\ell(\cdot)$ are differentiable.  From Equations (\ref{barepid-poisson}) and (\ref{bar-z}) it follows that
\begin{equation*}
|\bar Z_n(t)-\bar z(t)|\le |\frac{Z(0)}{n(t)}-\bar z_0e^{-(\lambda-\mu)t}| + \sum_\ell |\ell|\sup_{u\le t}|n^{-1}\hat Y\left(ne^{(\lambda-\mu)t}\bar\beta_\ell u \right)| + \int_0^tM|\bar Z_n(u)-\bar z(u)|du.
\end{equation*}
From Gronwall's inequality we then have
\begin{equation*}
|\bar Z_n(s)-\bar z(s)|\le \left( |\frac{Z(0)}{n(s)}-\bar z_0e^{-(\lambda-\mu)s}| + \sum_\ell |\ell|\sup_{u\le s}|n^{-1}\hat Y\left(ne^{(\lambda-\mu)s}\bar\beta_\ell u \right)| \right)e^{Ms}.
\end{equation*}
The exponent is independent of $n$ and both terms within the parenthesis tend to 0 with $n$ which implies that the left hand side tends to 0, also when taking supremum over finite intervals. This completes the proof.\hfill $\blacksquare$

\vskip.3cm
An alternative way of defining $\bar z(t)$ is in terms of differential equations rather than integral equations. If we differentiate Equation (\ref{bar-z}) we get, after a bit of algebra,
\begin{equation}
\bar z'(t)=\sum_\ell \ell \beta_\ell (\bar z(t)) - (\lambda-\mu)\bar z(t).
\end{equation}
Writing these term by term results in the following set of equations
\begin{align}
\bar s'(t)&=\lambda (1-\bar s(t))- \gamma \bar i(t)\bar s(t) \nonumber
\\
\bar i'(t)&=\bar i(t)\left( \gamma\bar s(t)-(\lambda+\delta)\right)\label{det-syst}
\\
\bar r'(t)&=-\lambda \bar r(t)+\delta\bar i(t). \nonumber
\end{align}

This set of differential equations are in fact well-known, they correspond to the SIR model with demography in a non-growing population (e.g.\ \cite[Section 4.4]{DHB12}. The starting point was defined as $\bar z_0=( 1-\varepsilon_1-\varepsilon_2, \varepsilon_1, \varepsilon_2)$, where $\varepsilon_1$ and $\varepsilon_2$ were assumed small but otherwise arbitrary. Because our main focus lies in the state of the process after a long time, the endemic level, we look for solutions $(\hat{\bar s}, \hat{\bar i}, \hat{\bar r})$ where all three derivatives are 0. Equating all derivatives to 0 gives the solutions $(1,0,0)$ and $(\hat{\bar s}, \hat{\bar i}, \hat{\bar r})$, where
\begin{equation}
\hat{\bar s}=\frac{\lambda+\delta}{\gamma},\qquad \hat{\bar i}=\lambda\left(\frac{1}{\lambda+\delta}-\frac{1}{\gamma}\right),\qquad \hat{\bar r}=\delta\left(\frac{1}{\lambda+\delta}-\frac{1}{\gamma}\right).\label{endlevel}
\end{equation}
Note that we have assumed that $\alpha> \lambda-\mu$ (super-critical case) ensuring that $0<\hat{\bar s}, \hat{\bar i}, \hat{\bar r}<1$, and we also have that $\hat{\bar s}+\hat{\bar i}+\hat{\bar r}=1$. Now, irrespective of the starting point, as long as $\bar i(0)>0‚ \ \bar s(0)\ge 0, \ \bar r(0)\ge 0$ and $\bar s(0)+\bar i(0)+\bar r(0)=1$, we have that the disease free solution $(1,0,0)$ is unstable, and therefore not a limit point of the process. For relevant parameter values (e.g.\ that the average infectious period is much shorter than the average life-length) this system  (\ref{det-syst}) is known to exhibit damped oscillations and the process converges to the endemic level: $(\bar s(t),\bar i(t),\bar r(t))\to (\hat{\bar s},\hat{\bar i},\hat{\bar r})$ as $t\to\infty$, \cite[Section 4.4]{DHB12}.

If we combine our results from the current and previous section for the case $\alpha > \lambda-\mu$ we have that $I_n(t)$ first behaves like $I^{(U)}(t)$. It may hence die out quickly, which it does with a probability tending to $(\delta+\mu)/\gamma$. With the remaining probability (tending to) $1-(\delta+\mu)/\gamma$, $I_n(t)$ grows at a higher rate than $N_n(t)$: $I_n(t)\sim e^{(\gamma-(\delta+\mu ))t}$ vs $N_n(t)\sim ne^{(\lambda-\mu)t}$. The sandwich coupling of $I_n(t)$ between $I^{(L)}(t)$ and $I^{(U)}(t)$ breaks down when the number of infected exceeds $\varepsilon_0 N_n(t)\sim \varepsilon_0 ne^{(\lambda-\mu)t}$ for the first time, which hence equals the first time that $I_n(t)= \varepsilon_0 N_n(t) \geq \varepsilon_0 ne^{(\lambda-\mu)t}$. From the growth rates just mentioned it can be deduced that this time $t_n$ satisfies 
\begin{equation}
t_n\approx \frac{\log (\varepsilon_0 n)}{\gamma-(\delta +\lambda )}.
\end{equation} 
(Note that the denominator is strictly positive by assumption.) At this time-point we have 
\begin{equation}
I_n(t_n)\approx (\varepsilon_0 n)^{(\gamma-(\delta+\mu))/(\gamma-(\delta +\lambda))}.
\end{equation} 
By our assumptions that $\lambda>\mu$ and $\gamma> \delta + \lambda$ it follows that the exponent is larger than 1. So, if $\varepsilon_1$ and $\varepsilon_2$ (the starting point of the approximation of the "proportion-processes" of infectives and recovered) are chosen small enough, then the second approximation kicks in before the first coupling is broken. As a consequence we then have convergence to endemicity. We summarize the results in the following theorem
 \begin{theorem}
Assume that $\alpha > \lambda-\mu$, let $Z_n(t)=(S_n(t), I_n(t),R_n(t))$ be the epidemic process starting with $(S_n(0),I_n(0),R_n(0))=(n-1,1,0)$ and  let $\bar Z_n(t)=Z_n(t)/ne^{(\lambda-\mu)t}$. Then, for fixed $n$ and as $t\to\infty$  we have that $Z_n(t)\to (\infty, 0, 0)$ (the disease free state in a growing population),  or else that $Z_n(t)\to (\infty,\infty,\infty)$. The probabiliy of the first event $\pi_n:= P(\lim_t Z_n(t)\to (\infty, 0, 0))$, satisfies $\lim_n\pi_n\to (\delta+\mu)/\gamma$. \end{theorem} 

In the latter event in the theorem we conjecture that, $\bar Z_n(t)\to (\hat{\bar s},\hat{\bar i},\hat{\bar r})$ as $t\to\infty$  (the endemic state in a growing population), this is true if the system (\ref{det-syst}) has no limit cycles in the positive octant. 

\emph{Remark}.
The system  (\ref{det-syst}) does not allow for $i(t) =0$ if the epidemic takes of. Since $n(t)$ is increasing  this implies that $I(t)$ can with large probability be kept above every desired minimum by choosing $n(0) = N_0(t)$ large enough. This implies that $I_n(t) \to \infty$ as $t\to \infty$ as claimed in part iv of Theorem \ref{mainthm}.

\section{Extension to the uniform Markovian dynamic SEIR epidemic model}\label{ext} 
Most infectious diseases have a latent state, i.e.\ a state where an individual has been infected (Exposed) but is not yet able to spread the disease, before becoming infectious. Accordingly we now extend our model to an SEIR epidemic allowing for this (the added ''E'' is for exposed). We have the same population model as before. The only difference is in the epidemic model that an individual who becomes infected is first Exposed (or latent) for an exponentially distributed time with rate $\nu$ (unless it dies before this) before turning infectious. After the latent period the individual becomes infectious having infectious contacts just like before. If we let $\nu\to\infty$ we hence retrieve the SIR model defined above. The transitions and their rates are given in Table \ref{Tab-SEIR}.

\begin{table}\centering
\begin{tabular}{|l|c|c|c|}
  \hline
Event & State change $\ell$& Rate name & Rate\\
  \hline
  Birth & $(1,0,0,0)$ & $\beta_{(1,00,0)}(s,i,r)$ & $\lambda (s+e+i+r)$\\
Death of susceptible & $(-1,0,0,0)$ &  $\beta_{(-1,0,0,0)}(s,i,r)$ & $\mu s$\\
Death of exposed & $(0,-1,0,0)$ & $\beta_{(0,-1,0,0)}(s,i,r)$ & $\mu e$\\
Death of infective & $(0,0,-1,0)$ & $\beta_{(0,0,-1,0)}(s,i,r)$ & $\mu i$\\
Death of recovered & $(0,0,0,-1)$ & $\beta_{(0,0,0,-1)}(s,i,r)$ & $\mu r$\\
Infection & $(-1,1,0,0)$ & $\beta_{(-1,1,0,0)}(s,i,r)$ & $\gamma i\cdot s/(s+e+i+r)$\\
Infectious & $(0,-1,1,0)$ & $\beta_{(0,-1,1,0)}(s,i,r)$ & $\nu e$\\
Recovery & $(0,0,-1,1)$ & $\beta_{(0,0,-1,1)}(s,i,r)$ & $\delta i$\\
  \hline
\end{tabular}
 \caption{The uniform Markovian dynamic SEIR epidemic: type of events, their state change $\ell$ (the old state $z=(s,e,i,r)$ is hence changed to $z-\ell$) and their rates.}\label{Tab-SEIR}
\end{table}

It is not hard to show that the basic reproduction number, defined as the expected number of infectious contacts an infected individual has, equals
\begin{equation}
R_0=\frac{\nu}{\nu+\mu}\frac{\gamma}{\delta+\mu}.\label{R-SEIR}
\end{equation}
The first factor is the probability that the individual does not die during the latent state, and the second is, as before, the expected number of infectious contacts while being infectious. As for the Malthusian parameter $\alpha$, the exponential rate at which the epidemic initially grows, it is obtained from the defining equation $\int_0^\infty e^{-\alpha t}c(t)dt=1$, where $c(t)$ is the expected rate at which an infected individual has infectious contacts $t$ time units after he/she was infected. In the SEIR model, this contact rate equals $\gamma$, but only if the latent period has ended and the infectious period has not ended yet (otherwise it equals 0). By conditioning on when the latent period ends, it follows that $c(t)=\gamma\nu\left( e^{-(\mu+\delta)t}-e^{-(\mu+\nu)t}\right) /(\nu-\delta)$. Inserting this into the defining equation for the Malthusian parameter $\alpha$ shows that it $\alpha $ is given by
\begin{equation}
\alpha = -\left( \mu +\frac{\delta+\mu}{2}\right) + \sqrt{ \frac{(\delta-\nu)^2}{4} + \gamma\nu}.\label{Malt-SEIR}
\end{equation}
It is easy to show that, as before, $R_0>1$ if and only if $\alpha >0$, and similarly if we replace ''$>$'' by ''='' or  ''$<$''.

It is possible to prove the same type of results as for the SIR epidemic. The same method of proofs apply. The only difference for the initial phase is that we should use the new $\alpha$ defined in (\ref{Malt-SEIR}), and that we get a new expression for the limit of $\pi_n$, the probability for not having a major outbreak, which now equals 
\begin{equation}
\lim_{n \to\infty}\pi_n=\frac{\delta+\mu}{\gamma} + \frac{\mu}{\nu+\mu}. \label{ext-SEIR}
\end{equation}
This expression is obtained by first deriving the offspring distribution which now is a mixture of a point mass at 0 (if the infected person dies during the latent state), and geometric as before. The extinction probability is derived using methods from branching process theory \cite{Jage75} and it is straightforward to show that it equals the right hand side of (\ref{ext-SEIR}).  For the SIR epidemic the minor outbreak probability was only the first term, which we retrieve if we make the latency period shrink down to 0 by letting $\nu\to\infty$.

For the endemic situation we now study the vector process $$\bar Z_n(t)=(S_n(t), E_n(t), I_n(t), R_n(t))/n(t)$$ for the same $n(t)=ne^{(\lambda-\mu)t}$ as before. The corresponding set of differential equations are
\begin{align}
\bar s'(t)&=\lambda (1-\bar s(t))- \gamma \bar i(t)\bar s(t) \nonumber
\\
\bar e'(t)&=\gamma \bar i(t)\bar s(t) -(\lambda+\nu)\bar e(t)\nonumber
\\
\bar i'(t)&=  \nu \bar e(t)-(\lambda +\delta)\bar i(t) \label{SEIR-det-syst}
\\
\bar r'(t)&=\delta\bar i(t) -\lambda \bar r(t). \nonumber
\end{align}
As in the SIR case, this set of differential equations correspond to the SEIR model with demography for a non-growing population \cite[Exercise 4.8]{DHB12}.
The endemic equilibrium (when $\alpha>\lambda-\mu$) is obtained by setting all derivatives equal to 0 and finding the unique positive solution. It is given by
\begin{equation}
(\hat{\bar s}, \hat{\bar e}, \hat{\bar i}, \hat{\bar r})=
\left( \frac{1}{\gamma b}, \
\frac{\lambda(\lambda+\delta)}{\nu}\left(b-\frac{1}{\gamma}\right), \
\lambda \left(b-\frac{1}{\gamma}\right), \
\delta\left(b-\frac{1}{\gamma}\right)\right) ,\label{SEIRendlevel}
\end{equation}
where $b=\nu/((\lambda+\nu)(\lambda+\delta))$. The theorems for the endemic situation apply also to the SEIR case using identical method of proofs, but with the new $\alpha$ from Equation (\ref{Malt-SEIR}), and for the new vector process $\bar Z(t)$ and its deterministic limit $\bar z(t)$ with its endemic equilibrium just defined. As for the initial phase the limiting probability of a minor outbreak $\lim_n\pi_n$ is different from the SIR model and equal to $(\delta+\mu)/\gamma + \mu/(\nu+\mu)$ as mentioned above.

\bibliographystyle{siam}
\bibliography{publications2}

\begin{thebibliography}{10}

\bibitem{AB00}
{\sc H.~Andersson and T.~Britton}, {\em Stochastic epidemic models and their
  statistical analysis}, vol.~4, Springer New York, 2000.

\bibitem{Ball95}
{\sc F.~Ball and P.~Donnelly}, {\em Strong approximations for epidemic models},
  Stochastic Process. Appl., 55 (1995), pp.~1--21.

\bibitem{BVB05}
{\sc R.~Breban, R.~Vardavas, and S.~Blower}, {\em Linking population-level
  models with growing networks: a class of epidemic models}, Physical Review E,
  72 (2005), p.~046110.

\bibitem{DHB12}
{\sc O.~Diekmann, H.~Heesterbeek, and T.~Britton}, {\em Mathematical Tools for
  Understanding Infectious Disease Dynamics}, Princeton University Press, 2012.

\bibitem{Durr07}
{\sc R.~Durrett}, {\em Random graph dynamics}, vol.~20, Cambridge university
  press, 2007.

\bibitem{EK05}
{\sc S.~N. Ethier and T.~G. Kurtz}, {\em Markov processes: characterization and
  convergence}, vol.~282, Wiley, 2005.

\bibitem{Hacc05}
{\sc P.~Haccou, P.~Jagers, and V.~Vatutin}, {\em {Branching processes:
  Variation, growth, and extinction of populations}}, Cambridge University
  Press, 2005.

\bibitem{Jage75}
{\sc P.~Jagers}, {\em Branching Processes with Biological Applications}, Wiley,
  New York, 1975.

\bibitem{Li99}
{\sc M.~Y. Li, J.~R. Graef, L.~Wang, and J.~Karsai}, {\em Global dynamics of a
  seir model with varying total population size}, Mathematical Biosciences, 160
  (1999), pp.~191--213.

\bibitem{Thie92}
{\sc H.~R. Thieme}, {\em Epidemic and demographic interaction in the spread of
  potentially fatal diseases in growing populations}, Mathematical biosciences,
  111 (1992), pp.~99--130.

\end{thebibliography}
\end{document}